\documentclass[preprint,12pt]{elsarticle}

\usepackage{babel}

\usepackage{hyperref}

\usepackage{graphicx}
\usepackage[left=3cm,right=1.5cm, top=2cm,bottom=2cm,bindingoffset=0cm]{geometry}
\usepackage[cp1251]{inputenc}
\usepackage{amssymb}
\usepackage{amsthm}
\usepackage{cmap}
\usepackage{epsfig}
\usepackage{amsmath}

\biboptions{sort&compress}

\newcommand{\sysn}{\left\{\begin{array}{rcl}}
\newcommand{\sysk}{\end{array}\right.}

\renewcommand{\le}{\leqslant}

\newtheorem{theorem}{Theorem}[section]

\newtheorem{example}[theorem]{Example}

\newtheorem{definition}[theorem]{Definition}
\newtheorem{remark}[theorem]{Remark}

\journal{...}

\begin{document}

\title{Various $S(n)$-closednesses in $S(n)$-spaces with examples}

\author{Alexander V. Osipov}
\address{Krasovskii Institute of Mathematics and Mechanics, Ural Federal
 University, Yekaterinburg, Russia}

\ead{OAB@list.ru}

\begin{abstract}  In this paper we continue to study various types of closures in $S(n)$-spaces. The main results are related to the construction and illustration of examples that allow us to understand the relationship between $S(n)$-closed, $S(n)$-$\theta$-closed, weakly $S(n)$-closed and weakly $S(n)$-$\theta$-closed spaces for each $n\in \mathbb{N}$. The relation of these classes in Lindel\"{o}f spaces is shown. Some of the solved problems formulated by D. Dikranjan and E. Giuli are presented in the examples.
\end{abstract}


\maketitle

\section{Introduction}

In 1924, P.S. Alexandroff and P.S. Urysohn \cite{alur} established a number of characterizations of compactness including the following: { \it a space compact if and only if every infinite subset has a complete accumulation point}. Also in \cite{alur}, it introduced and characterized the concept of $H$-closed spaces. A Hausdorff space is said to be {\it $H$-closed} (or absolutely closed) if it is closed in every Hausdorff space containing it as a subspace. This property is a generalization of compactness, since a compact subset of a Hausdorff space is closed. Thus, every compact Hausdorff space is $H$-closed.

Alexandroff and Urysohn \cite{alur} extended the complete accumulation point characterization of compactness to $H$-closure by proving that any
$H$-closed space has the following property:

$(*)$ {\it any infinite set of regular power has a $\theta$-complete
accumulation point}, i.e.,  there is a point $p\in X$ such that for each neighborhood $U$ of $p$, $|A\cap \overline{U}|=|A|$.

However, the converse is not true. The first example of a space
possessing property $(*)$ and not being $H$-closed was constructed
by G.A. Kirtadze \cite{kir}.

\begin{example}(Example 3 in \cite{kir}) Let $T=((\omega_1 +1)\times (\omega +1)\setminus \{\omega_1,\omega \})$  be the deleted Tychonoff plane and let $X=T\times(\omega+1)$, whose elements will
be denoted by  $(\alpha,k,i)$ where $1\leq \alpha\leq \omega_1$,  $1\leq k\leq \omega$ and $1\leq i\leq \omega$.

Consider the following identifications on the space $X$ (see Figure~1):

1. The set of points of the form $(\alpha, k,\omega)$ with the third coordinate $\omega$ determines the point $\xi$. (In the drawing, this is the upper face of the 'parallelepiped').

2. For each natural number $n$, a countable set of points $(\omega_1, 1, n)$, $(\omega_1, 3, n)$, ..., $(\omega_1, 2k-1, n)$, ...
defines a point $\eta_n$.

3. For each natural number $n$, a finite set of points
$(\omega_1, 2n,1)$, $(\omega_1, 2n,2)$, ..., $(\omega_1, 2n, n)$
defines a point $\zeta_n$.

\end{example}

\begin{figure}[h!]
\centering
\includegraphics[width=0.7\textwidth]{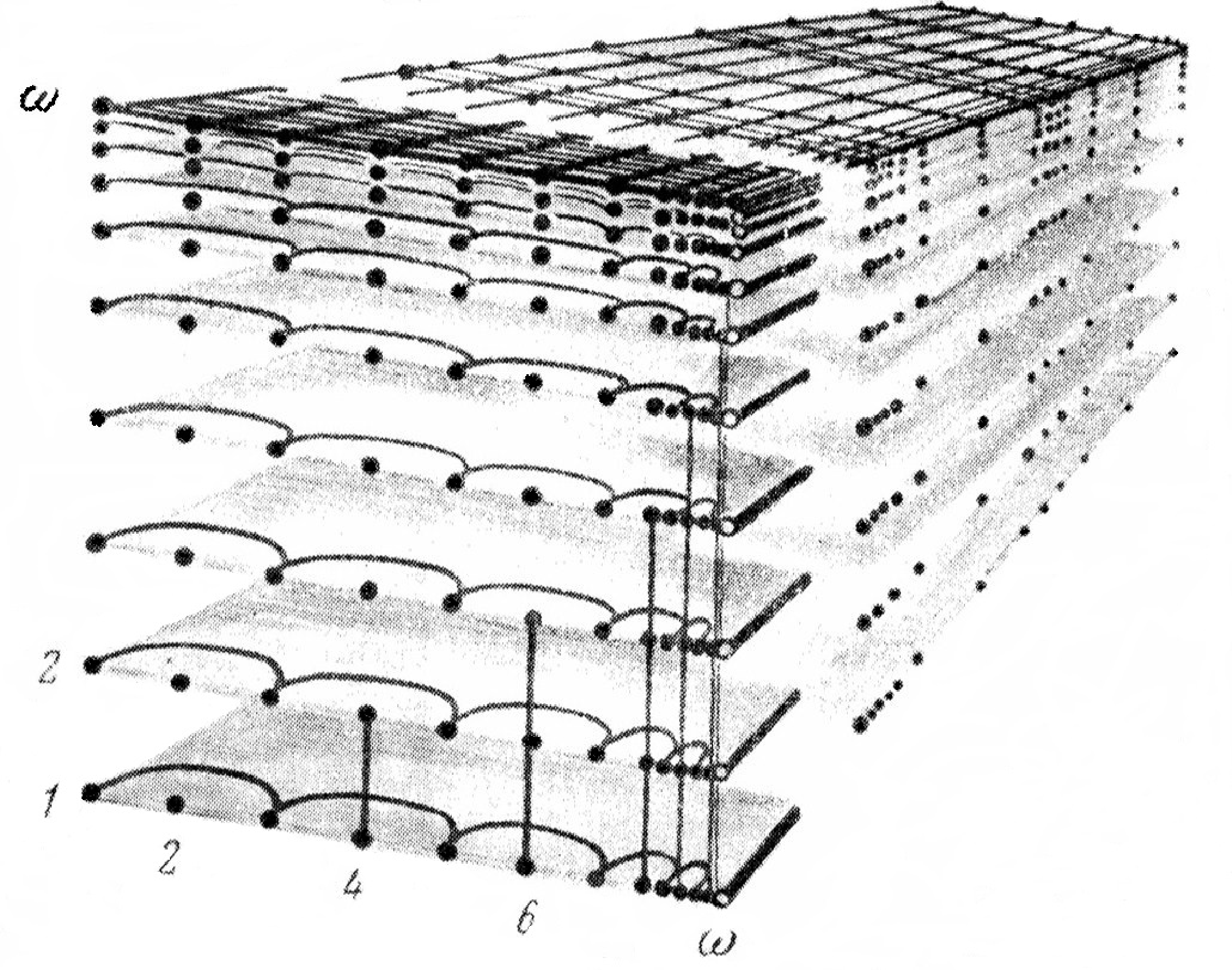} 
\caption{Kirtadze's space $K$}
\end{figure}



Let $K$ be the space $X=(\bigcup\limits_{i=1}^{\infty} T_i)\cup \{\xi\}$ with respect to these identifications where the point $\xi$ with the base of neighborhoods $U_n(\xi)=X\setminus \bigcup\limits_{i=1}^n  T_i$.

The obtained Hausdorff space $K$ has the property $(*)$ and is not $H$-closed (see \cite{kir}).

\bigskip

In \cite{osip1}, it constructed a simple example with these properties (see below Example \ref{1}).  Such simple example with these properties are also constructed in \cite{por1}.

\bigskip

In 1966, Velichko \cite{vel} introduced the notion of $\theta$-closedness. For a subset $M$ of a topological space $X$,
the $\theta$-closure $cl_{\theta} M$ is defined as the set of $x\in X$ such that any closed neighborhood of $x$ intersects $M$.
This notion has been used extensively to study non-regular Hausdorff spaces.

The $S(n)$-spaces were introduced by Viglino in 1969 (see \cite{vig}) under the name $\overline{T}_n$-spaces. After that $S(n)$-spaces and $S(n)$-closed spaces were studied by other authors. For example, J. Porter in 1969 (see \cite{port1}) studied minimal $R(\omega_0)$ spaces, where he used the notation $R(n)$ for $S(2n-1)$-spaces and $U(n)$ for $S(2n)$-spaces. For the first time the notation $S(n)$ for $S(n)$-spaces appeared in 1973 in \cite{porvot} where the authors extended the definition of $S(n)$-spaces to $S(\alpha)$-spaces, where $\alpha$ is any
ordinal. In that paper Porter and Votaw, among other results, characterized
the minimal $S(\alpha)$ and $S(\alpha)$-closed spaces. In 1986, Dikranjan and Giuli introduced the more general notion of $\theta^n$-closure and developed the theory of  $S(n)$-$\theta$-closed spaces \cite{dg}.
In 2003, it introduced the notions of weakly $S(n)$-closed and weakly $S(n)$-$\theta$-closed spaces and continued to study the theory of  $S(n)$-spaces~\cite{osip5}.
\bigskip

In this paper we continue to study various types of closures in $S(n)$-spaces. The main results in Section 3 are related to the construction and illustration of examples that allow us to understand the relationship between $S(n)$-closed, $S(n)$-$\theta$-closed, weakly $S(n)$-closed and weakly $S(n)$-$\theta$-closed spaces for each $n\in \mathbb{N}$. In Section 4, the relation of these classes in Lindel\"{o}f spaces is shown. Some of the solved problems formulated by Dikranjan and Giuli are presented in examples in Section 5.

Throughout the paper a space means a Hausdorff space; $\overline{M}$ or $cl M$ denotes the closure of the set $M$ in a given topological space; $\mathbb{N}$ denotes the set of positive integers. Recall that a subset $W$ of $X$ is called {\it functionally open (co-zero)} if there is a continuous function $f: X\rightarrow \mathbb{R}$ such that $W=f^{-1}(\mathbb{R}\setminus \{0\})$, i.e., $W$ is a complement of the zero-set $f^{-1}(0)$.

\section{Main definitions}

\begin{definition}(see \cite{dg}) Suppose that $X$ is a topological space, $M\subset X$, and $x\in X$. For each $n\in \mathbb{N}$, the {\it $\theta^n$-closure operator} is defined as follows: $x\notin cl_{\theta^n}M$ if there exists a set of open neighborhoods $U_1\subset U_2\subset...\subset U_n$ of the point $x$ such that $cl U_i\subset U_{i+1}$ for $i=1,2,...,n-1$ and $clU_n\cap M=\emptyset.$ For $n=0$, we put $cl_{\theta^0}=clM$.
\end{definition}

For $n=1$, this definition gives the $\theta$-closure operator defined by Velichko.

A set $M$ is said to be $\theta^n$-closed if $M=cl_{\theta^n} M$. Denote by $Int_{\theta^n} M=X\setminus cl_{\theta^n}(X\setminus M)$ the $\theta^n$-interior of the set $M$. Evidently, $cl_{\theta^n}(cl_{\theta^s} M)\subset cl_{\theta^{n+s}}M$ for $M\subset X$ and $n,s\in \mathbb{N}$. For $n\in \mathbb{N}$ and a filter $\mathcal{F}$ on $X$, denote by $ad_{\theta^n} \mathcal{F}$ the set of $\theta^n$-adherent points, i.e.,
$ad_{\theta^n} \mathcal{F}=\{\bigcap\limits_{\alpha} cl_{\theta^n} F_{\alpha}: F_{\alpha}\in \mathcal{F}\}$. In particular, $ad_{\theta^0} \mathcal{F}=ad \mathcal{F}$ is the set of adherent points of the filter the $\mathcal{F}.$ For any $n\in \mathbb{N}$, a point $x\in X$ is {\it $S(n)$-separated} from a subset $M$ if $x\notin cl_{\theta^n} M$. For example, $x$ is $S(0)$-separated from $M$ if $x\notin \overline{M}$. For $n>0$, the relation of $S(n)$-separability of points is symmetric. On the other hand, $S(0)$-separability may be not symmetric in some not $T_1$-spaces. Therefore, we say that points $x$ and $y$ are $S(0)$-separated if
$x\notin \{\overline{y}\}$ and $y\notin \{\overline{x}\}$.
Let $n\in \mathbb{N}$ and $X$ be a topological space.

\medskip
1. $X$ is called an {\it $S(n)$-space} if any two distinct points of $X$ are $S(n)$-separated.

\medskip

2. A filter $\mathcal{F}$ on $X$ is called an {\it $S(n)$-filter} if every point, not being an adherent point of the
filter $\mathcal{F}$, is $S(n)$-separated from some element of the filter $\mathcal{F}$.

\medskip

3. An open cover $\{U_{\alpha}\}$ of the space $X$ is called an {\it $S(n)$-cover} if every point of $X$ lies in the
$\theta^n$-interior of some $U_{\alpha}$.

\medskip

It is obvious that $S(0)$-spaces are $T_0$-spaces, $S(1)$-spaces are Hausdorff spaces, and $S(2)$-spaces
are Urysohn spaces. It is clear that every filter is an $S(0)$-filter, every open cover is an $S(0)$-cover,
and every open filter is an $S(1)$-filter. Open $S(2)$-filters are called Urysohn filters. $S(1)$-covers are called Urysohn covers. In a regular space,
every filter (every cover) is an $S(n)$-filter ($S(n)$-cover) for any $n\in \mathbb{N}$.

\begin{definition} Let $n\in \mathbb{N}$. A neighborhood $U$ of a point $x$ is called  an {\it $n$-hull} of the point $x$ if there exists
a set of neighborhoods $U_1$, $U_2$, ..., $U_n = U$ of the point $x$ such that $\overline{U_i}\subset U_{i+1}$ for $i=1,...,n-1$. In particular, an $1$-hull of the point $x$ is a neighborhood of $x$.
\end{definition}

\begin{definition} A point $x$ from $X$ is called

$\bullet$ {\it a $\theta^0(n)$-complete accumulation} point of
an infinite set $F$ if $|F\cap U|=|F|$  for arbitrary $n$-hull of the point $x$;

$\bullet$ {\it a  $\theta(n)$-complete accumulation} point of
an infinite set $F$ if $|F\cap \overline{U}|=|F|$ for arbitrary $n$-hull of the point $x$.
\end{definition}

\section{Various $S(n)$-closednesses in $S(n)$-spaces}

A topological $S(n)$-space $X$ is called

\medskip

$\bullet$ {\it $S(n)$-closed},  if it is closed in every $S(n)$-space containing it as a subspace;

\medskip

$\bullet$ {\it $S(n)$-$\theta$-closed},  if it is $\theta$-closed in every $S(n)$-space containing it as a subspace;

\medskip

$\bullet$  {\it weakly $S(n)$-$\theta$-closed}, if any infinite set of regular power of the space $X$ has a $\theta^0(n)$-complete accumulation point;

\medskip

$\bullet$  {\it weakly $S(n)$-closed}, if any infinite set of regular power of the space $X$ has a $\theta(n)$-complete accumulation point.

\medskip

Note that a $\theta^0(1)$-complete accumulation is a point of complete accumulation, and a $\theta(1)$-complete accumulation is a $\theta$-complete accumulation point.
Thus,  weakly $S(1)$-$\theta$-closed and weakly $S(1)$-closed spaces are compact Hausdorff spaces
and spaces with property $(*)$, respectively. Spaces with property $(*)$ (= weakly $S(1)$-closed spaces) is called {\it weakly $H$-closed} spaces \cite{osip6} (or {\it nearly $H$-closed} spaces \cite{osip16}).

\medskip
Note that every compact Hausdorff space is  $S(1)$-$\theta$-closed, hence the properties of $S(1)$-$\theta$-closed and weakly $S(1)$-$\theta$-closed are equivalent and equal to compactness.

\medskip
Porter and Votaw \cite{porvot} characterized $S(n)$-closed spaces by means of open $S(n)$-filters and $S(n)$-covers
(for $n=2$, see Herrlich \cite{herr}).

\medskip

Let $n\in \mathbb{N}^+$ and $X$ be an $S(n)$-space. Then the following conditions are equivalent:

(1) $ad_{\theta^n} \mathcal{F}\neq \emptyset$ for any open filter $\mathcal{F}$ on $X$;

(2) $ad \mathcal{F}\neq \emptyset$ for any open $S(n)$-filter $\mathcal{G}$ on $X$;

(3) for any $S(n-1)$-cover $\{U_{\alpha}\}$ of $X$ there exist $\alpha_1$,..., $\alpha_k$ such that $X=\bigcup\limits_{i=1}^k \overline{U}_{\alpha_i}$;

(4) $X$ is an $S(n)$-closed space.

\bigskip

Dikranjan and Giuli \cite{dg} characterized $S(n)$-$\theta$-closed spaces in terms of $S(n-1)$-filters and
$S(n-1)$-covers.

Let $n\in \mathbb{N}^+$ and $X$ be an $S(n)$-space. Then the following conditions are equivalent:

(1) $ad_{\theta^{n-1}} \mathcal{F}\neq \emptyset$ for any closed filter $\mathcal{F}$ on $X$;

(2) $ad \mathcal{F}\neq \emptyset$ for any closed $S(n-1)$-filter $\mathcal{G}$ on $X$;

(3) for any $S(n-1)$-cover of the space $X$ has a finite subcover;

(4) $X$ is an $S(n)$-$\theta$-closed space.

\medskip

Note that $S(1)$-closedness and $S(1)$-$\theta$-closedness are $H$-closedness and compactness,
respectively. $S(2)$-closedness and $S(2)$-$\theta$-closedness are $U$-closedness and $U$-$\theta$-closedness,
respectively. From characteristics themselves, it follows that

\medskip

${\bf (P0)}$ a $S(n)$-$\theta$-closed subspace of an $S(n)$-
space is an $S(n)$-closed space.

\medskip

The following assertions are holds for every $n>1$.

\medskip

${\bf (P1)}$  $S(n-1)$-closedness yields $S(n)$-$\theta$-closedness (Corollary 2.3. in ~\cite{dg}).

\medskip

${\bf (P2)}$   $S(n)$-$\theta$-closedness yields weakly $S(n)$-$\theta$-closedness (Theorem 1 in ~\cite{osip5}).

\medskip

${\bf (P3)}$   $S(n)$-closedness yields weakly $S(n)$-closedness (Theorem 2 in ~\cite{osip5}).

\medskip
Note that any $\theta^0(n)$-complete accumulation point is a $\theta(n)$-complete accumulation point, it follows that

\medskip

${\bf (P4)}$  weakly $S(n)$-$\theta$-closedness yields weakly $S(n)$-closedness.

\medskip

 Moreover, since a $\theta(n)$-complete accumulation point is a $\theta^0(n+1)$-complete accumulation point, it follows that

 \medskip

${\bf (P5)}$  weakly $S(n)$-closedness yields weakly $S(n+1)$-$\theta$-closedness.

\bigskip

Below we present examples that separate the studied classes of $S(n)$-spaces.

 \begin{example}\label{1} (Example 1 in \cite{osip7, osip5}) Let $T_1$ and $T_2$ be two copies of the deleted Tychonoff
 plane  $T$, whose elements will be denoted by $(\alpha,n,1)$ and
$(\alpha,n,2)$, respectively. On the topological sum $T_1\oplus T_2$, we consider the
identifications

\medskip

    $(\omega_1,k,1)\sim (\omega_1,2k,2)$ for every $k\in
    \mathbb{N}$;

    and we identify all points  $(\omega_1,2k-1,2)$ for any  $k\in \mathbb{N}$ with the same
    point ${\bf b}$.

    Adding, to the obtained space, a point ${\bf a}$ with the base of neighborhoods

 $U_{\beta,k}({\bf a})=\{(\alpha,n,1): \beta < \alpha <\omega_1,k<n\le \omega_0 \}\cup\{{\bf a}\}$
 for  arbitrary $\beta < \omega_1$  and $ k<\omega_0 $, we get a
 space $X_1$ (see Figure~2) with the following properties.
\end{example}
\bigskip

\begin{figure}[h!]
\centering
\includegraphics[width=0.8\textwidth]{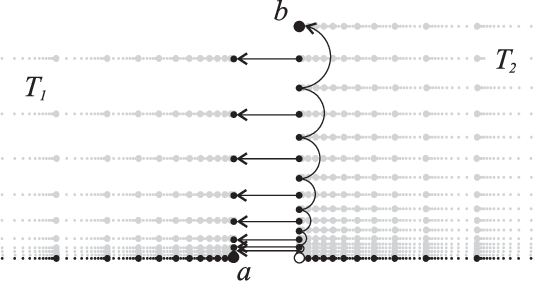} 
\caption{The space $X_1$}
\end{figure}




\bigskip

 $\bullet$ The space $X_1$ is an Urysohn space, i.e., $S(2)$-space.

 \medskip

 $\bullet$ Note that the closure of any neighborhood $U_{\beta,k}({\bf a}) $ of ${\bf a}$
contains all points $(\omega_1,n,1)$  but a finite number, and
the closure of any neighborhood of the point ${\bf b}$ contains
all points $(\alpha,\omega,2)$  but a countable number. This
evidently implies that $X_1$ is a weakly $H$-closed space.

\medskip

$\bullet$  Considering a system of open sets with the finite
intersection property
       $V_{\alpha,k}=\{(\beta,p,2): \alpha < \beta <\omega_1,$
    for all even $p>k \} $ for arbitrary $\alpha$ and $k$, we
    obtain  $\bigcap\limits_{\alpha,\,k}\overline{V_{\alpha,\,k}}=\emptyset$.
      Thus, {\it the Urysohn space $X_1$ is an example of a space not $H$-closed but weakly
      $H$-closed.}

 \bigskip

\begin{remark} Identifying the points ${\bf a}$ and ${\bf b}$ in
$X_1$, we obtain a space $X_1^1$ with the following properties.

$\bullet$ $X_1^1$ is functionally Hausdorff, i.e.,  any two distinct points can be separated by a continuous function (there exists a continuous function $f: X\rightarrow[0,1]$ with $f(x)=0$ and $f(y)=1$). This follows from the fact that for any two points $x$ and $y$ of the space $X_1^1$ there exists an open-closed subset $W$ such that $x\in W$ and $y\notin W$. Then $f:X\rightarrow[0,1]$ such that $f(W)=0$ and $f(X\setminus W)=1$ is the required function.

$\bullet$  $X_1^1$ is a $CH$-closed space, i.e., a functionally Hausdorff space such that it is closed in every functionally Hausdorff space in which it can be embedded.

Indeed, the complement of any functionally open neighborhood of $c=\{a,b\}$ is compact. This means that any cover of $X_1^1$  by functionally open sets has a finite subcover. This is equivalent to being $CH$-closed (Theorem 4.9 in  \cite{bps}). Note that being $CH$-closed is also equivalent to the Stone-Weierstrass theorem being satisfied on the space (Theorem 4.9 in \cite{bps}).

$\bullet$  $X_1^1$ is non-$H$-closed.  Consider the system of open sets with the finite intersection property  $V_{\alpha,k}=\{(\beta,p,2): \alpha < \beta <\omega_1,$   for all even $p>k \} $ for arbitrary $\alpha$ and $k$,
 we obtain $\bigcap\limits_{\alpha,\,k}\overline{V_{\alpha,\,k}}=\emptyset$.

Thus, the space $X_1^1$ is an example of not $H$-closed, but $CH$-closed space on which the Stone-Weierstrass theorem holds.

\end{remark}

\bigskip
\begin{example}(Example 2 in \cite{osip5}) Let $n>1$ and let $T^i_1$ and  $T^i_2$ ( $i=1,...,n$ ) be  $2n$
copies of the deleted Tychonoff plane $T$, whose elements will be denoted
by $(i,\alpha,k,1)$ and $(i,\alpha,k,2)$, respectively. Consider
the following identifications on the topological sum  $(\bigoplus\limits_{i=1}^n T^i_1) \bigoplus (\bigoplus\limits_{i=1}^n T^i_2):$

  $(1,\omega_1,k,1)\sim (1,\omega_1,2k,2)$ for every $k\in \mathbb{N}$;

  $(s,\alpha,\omega_0,1)\sim (s+1,\alpha,\omega_0,1)$ for odd
  $s$;

  $(s,\omega_1,k,1)\sim (s+1,\omega_1,k,1)$ for even $s$;

  $(1,\omega_1,2k-1,2)\sim (2,\omega_1,k,2)$ for every $k\in \mathbb{N}$;

  $(s,\alpha,\omega_0,2)\sim (s+1,\alpha,\omega_0,2)$ for even
  $s$;

  $(s,\omega_1,k,2)\sim (s+1,\omega_1,k,2)$ for odd $s>1$.

Adding, to the space obtained, two points ${\bf a}$ and ${\bf b}$
with the base of neighborhoods:

 $U_{\alpha,k}({\bf a})=\{(n,\beta,s,1): \alpha < \beta <\omega_1, k<s<\omega_0
 \}\cup\{{\bf a}\}$, and

 $U_{\alpha,k}({\bf b})=\{(n,\beta,s,2): \alpha < \beta <\omega_1, k<s<\omega_0 \}\cup\{{\bf b}\}$

 for any $\alpha<\omega_1$  and $k<\omega_0$, we get an
 $S(n)$-space $X_n$ (see Figure~3 for $n=4$), which will be weakly $S(n)$-closed.
\end{example}

\bigskip

\begin{figure}[h!]
\centering
\includegraphics[width=1\textwidth]{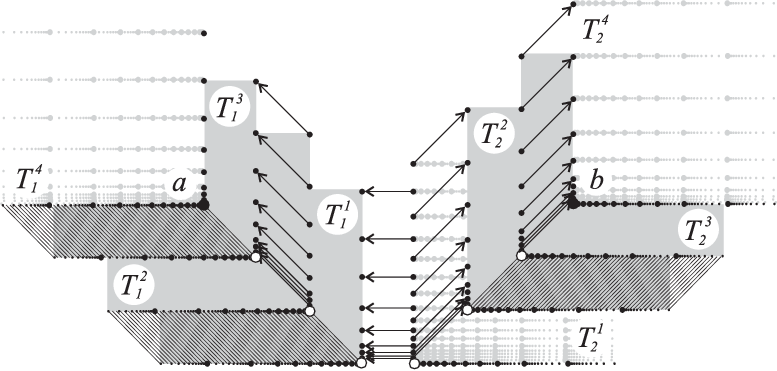} 
\caption{The space $X_4$}
\end{figure}



Indeed, the closure of any $n$-hull of the point ${\bf a}$
contains all the sets $\{(i,\beta,\omega_0,1): i=1,...,n$ and
$\beta<\omega_1\}$  but a countable  number of points, and all the sets
$\{(i,\omega_1,k,1): i=1,...,n$ and $k<\omega_0\}$  a finite number of points, and the closure of an $n$-hull of the point ${\bf
b}$ contains all the sets $\{(i,\beta,\omega_0,2): i=1,...,n$ and
$\beta<\omega_1\}$ but a countable number of points,  and all the sets
$\{(i,\omega_1,k,2): i=1,...,n$ and $k<\omega_0\}$ but a finite number of points and except the points $(1,\omega_1,2k,2)$.

 Considering an open filter with the base
 $V_{\alpha, k}=\{(1,\beta,p,2): \alpha < \beta <\omega_1$ and for
all even $ p>k \} $ for arbitrary $\alpha$ and $k$, we obtain
$ad_{\theta^n}\{V_{\alpha,k}\}_{\alpha,k}=\emptyset$.

Hence, the $S(n)$-space $X_n$ is weakly $S(n)$-closed but not
$S(n)$-closed.

\begin{example}\label{124} (Example 3 in \cite{osip5}) Let $n>1$ and let $T^i_1$, $T_2$, $T^i_3$ ( $i=1,...,n-1$ )
be   $2n-1$ copies of the deleted Tychonoff plane $T$, whose elements will
be denoted by $(i,\alpha,k,1)$, $(\alpha,k,2)$ and
 $(i,\alpha,k,3)$, respectively.

Consider the following identifications on the topological sum $(\bigoplus\limits_{i=1}^{n-1} T^i_1) \bigoplus T_2 \bigoplus (\bigoplus\limits_{i=1}^{n-1} T^i_3):$

$(1,\omega_1,k,1)\sim (\omega_1,k,2)$ for every $k\in \mathbb{N}$;

$(1,\alpha,\omega_0,3)\sim (\alpha,\omega_0,2)$ for every
$\alpha$;

$(s,\alpha,\omega_0,1)\sim (s+1,\alpha,\omega_0,1)$ for all odd
$s$ and every $\alpha$;

$(s,\omega_1,k,1)\sim (s+1,\omega_1,k,1)$ for all even $s$ and
every $k$;

 $(s,\omega_1,k,3)\sim (s+1,\omega_1,k,3)$ for all odd $s$ and
 every $k$;

 $(s,\alpha,\omega_0,3)\sim (s+1,\alpha,\omega_0,3)$ for all even $s$ and
every $\alpha$.

To the space obtained, we add two points ${\bf a}$ and ${\bf b}$
with the base of neighborhoods:

 $U_{\alpha,k}({\bf a})=\{(n-1,\beta,p,1):
\alpha < \beta <\omega_1, k<p<\omega_0 \}\cup\{{\bf a}\}$

and $U_{\alpha,k}({\bf b})=\{(n-1,\beta,p,3):
\alpha < \beta <\omega_1, k<p<\omega_0 \}\cup\{{\bf b}\} $ for
arbitrary $\alpha < \omega_1$  and $k<\omega_0$.

 We change the base of neighborhoods of the points $(\omega_1,k,2)$ and $(\alpha,\omega_0,2)$ for all $\alpha$
and $k$, setting  $U_{\alpha,k}((\omega_1,k,2))=\{\omega_1,k,2\}\bigcup V_{\alpha,k}\setminus\{(\alpha,k,2):
(\alpha,k,2)$ is not an isolated point in $T_2 \}$, where $V_{\alpha,k}$ is a standard neighborhood in the quotient
topology.
\bigskip

\end{example}

The obtained space $Y_n$ (see Figure~4 for $n=4$) has the following
properties:

\bigskip

\begin{figure}[h!]
\centering
\includegraphics[width=1.05\textwidth]{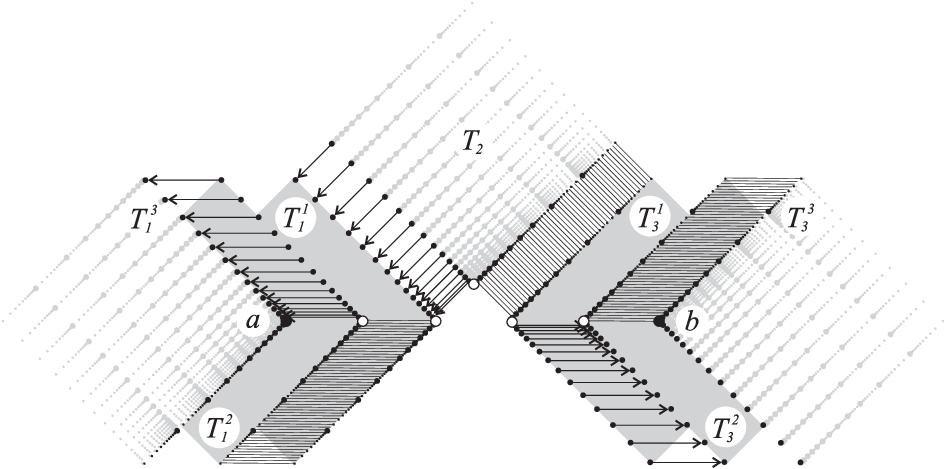} 
\caption{The space $Y_4$}
\end{figure}



 $\bullet$    $Y_n$  is an $S(n)$-space;

 \medskip

$\bullet$ the closure of an $n$-hull of the point ${\bf a}$
contains all points $(s,\alpha,\omega_0,1)$ for $s=1,...,n-1$  but a countable number, and all points $(s,\omega_1,k,1)$ for $s=1,...,n-1$,
but a finite number.

$\bullet$ the closure of an $n$-hull of the point ${\bf b}$
contains all points $(s,\alpha,\omega_0,3)$ for $s=1,...,n-1$ but a countable number, and all points $(s,\omega_1,k,3)$ for $s=1,...,n-1$,
but a finite number.

\medskip

Thus, the space $Y_n$ is weakly $S(n)$-$\theta$-closed.

Considering the base $\{F_{\alpha,k}\}$ of a closed filter ${\cal
F}$, where   $F_{\alpha, k}=\{(\beta,p,2): \alpha < \beta
<\omega_1$ and $(\beta,p,2)$ is not an isolated point in $T_2 \}$,
we get $ ad_{\theta^{n-1}} \cal F=\emptyset$. Hence, $Y_n$ is not
an $S(n)$-$\theta$-closed space.

Note that from the fact that any set lying in $T_2\setminus  P$,
where $P=\{(\alpha,k,2):\alpha=\omega_1$ or $k=\omega_0 \}$, has a
$\theta(1)$-complete accumulation point (lying in $P$) and any
other set (of regular power) has a $\theta(n-1)$-complete
accumulation point (either the point ${\bf a}$ or the point ${\bf
b}$), it follows that $Y_n$ is a weakly $S(n-1)$-closed space.

Thus, $Y_n$ is an example of an $S(n)$-space,
which, while being a weakly $S(n)$-$\theta$-closed space (even, moreover, $Y_n$ is weakly $S(n-1)$-closed space), is not an $S(n)$-$\theta$-closed space.

\begin{example}\label{125} (Example 4 in \cite{osip5}) Let $n>1$ and $T^i$ ( $i=1,...,n-1$ ) be  $n-1$
copies of the deleted Tychonoff plane $T$, whose elements will be denoted
$(i,\alpha,k)$. Let $\omega_1$ be the set of all countable
ordinals with the order topology, whose elements will be denoted
by  $\{\alpha\}$.

Consider the following identifications on the topological sum $\bigoplus\limits_{i=1}^{n-1}T^i\bigoplus\omega_1$:

\medskip

  $(1,\alpha,\omega_0)\sim (\{\alpha\})$, where $\alpha$ is the limit ordinal number in  $\omega_1$;

\medskip

 $(s,\omega_1,k)\sim (s+1,\omega_1,k)$ for all odd $s$;

\medskip

 $(s,\alpha,\omega_0)\sim (s+1,\alpha,\omega_0)$ for all even  $s$.

To the obtained space, we add the point ${\bf a}$ with the base of
neighborhoods:

 $U_{\alpha,k}({\bf a})=\{(n-1,\beta,p): \alpha < \beta <\omega_1,k<p<\omega_0
 \}\cup\{{\bf a}\}$.

\end{example}

\begin{figure}[h!]
\centering
\includegraphics[width=1\textwidth]{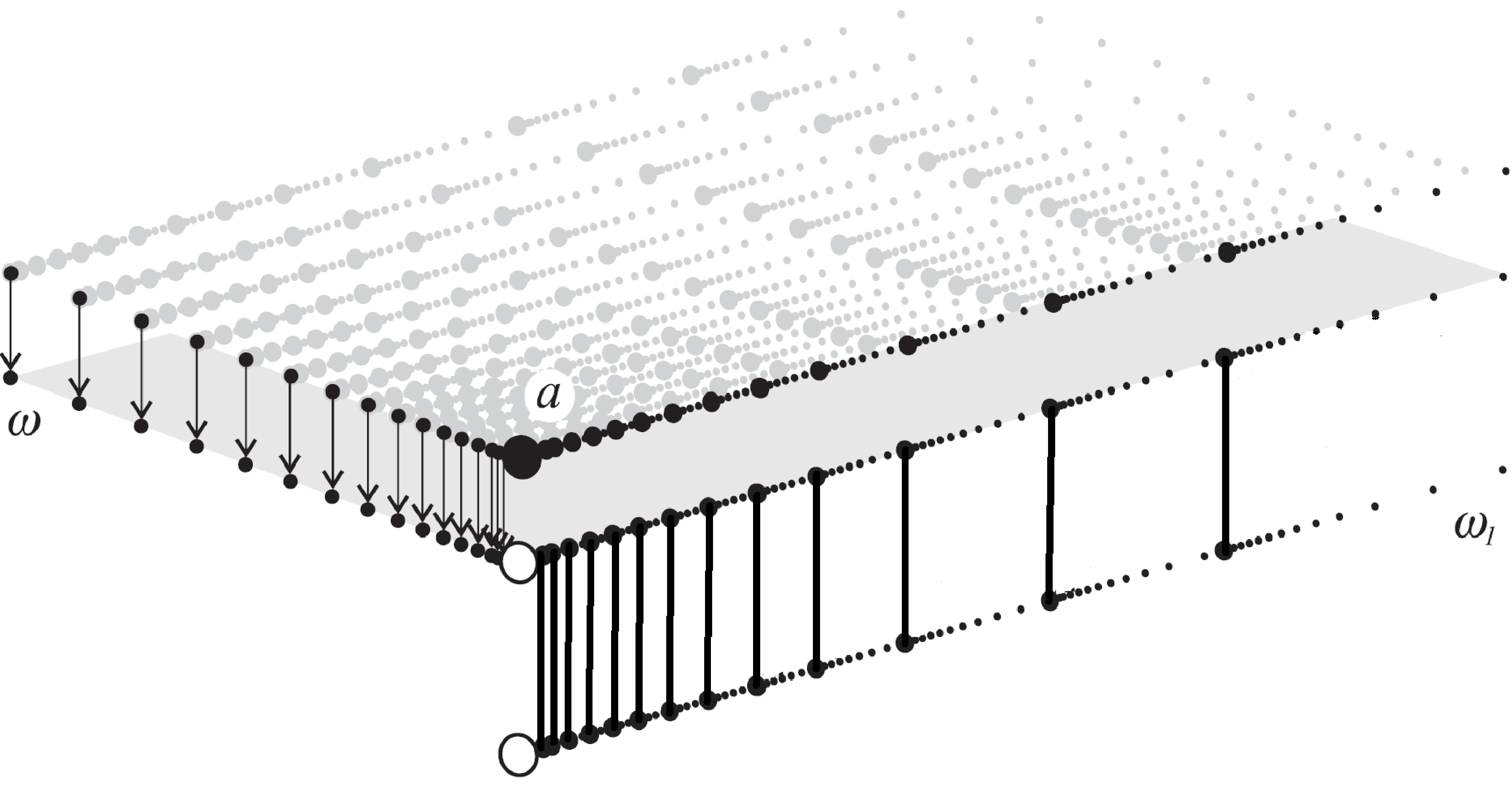} 
\caption{The space $Z_3$}
\end{figure}



The space $Z_n$ (Figure~5 for $n=3$) thus obtained has the following
properties:

$\bullet$ $Z_n$ is $S(k)$-space for any $k\in \mathbb{N};$

$\bullet$ $Z_n$ --- $S(n)$-$\theta$-closed space.

  Indeed, consider any $S(n-1)$-cover $V=\{V_{\alpha}\}$. Then ${\bf a}\in Int_{\theta^{n-1}} V_{\alpha}$ for some $V_{\alpha}$ from $V$.
 Hence, $V_{\alpha}$ is an $n$-hull of the point ${\bf a}$. However,
 the complement $\omega_1\setminus V_{\alpha}$ of the $n$-hull of the point ${\bf a}$ contains no more than countable set and, evidently, is covered by
 a finite number of elements of the $S(n-1)$-cover $V$. On the
 other hand, the set $F=\{\gamma: \gamma$ is isolated in $\omega_1
  \}$ is of the regular power $\omega_1$ and has no $\theta(n-1)$-complete accumulation point, since the closure of an
  $(n-1)$-hull of the point ${\bf a}$ does not intersect  $F$.
  Thus, the space $Z_n$ is $S(n)$-$\theta$-closed but not weakly $S(n-1)$-closed.

\bigskip

\begin{example}\label{126}  Let $n>1$ and let $T^i$ ($i=1,...,n$) be  $n$ copies of the deleted Tychonoff plane $T$, whose elements will
be denoted by  $(i,\alpha,k)$.

Consider the following identifications on the topological sum  $\bigoplus\limits_{i=1}^{n}T^i$:

 $(s,\omega_1,k)\sim (s+1,\omega_1,k)$ for all odd $s$;

  $(s,\alpha,\omega_0)\sim (s+1,\alpha,\omega_0)$ for all even  $s$.

To the obtained space, we add the point ${\bf a}$ with the base of
neighborhoods:

 $U_{\alpha,k}({\bf a})=\{(n,\beta,p): \alpha < \beta <\omega_1,k<p<\omega_0
 \}\cup\{{\bf a}\}$.

\bigskip
\end{example}

The space $G_n$ (see Figure~6 for $n=4$) thus obtained has the following
properties:

\begin{figure}[h!]
\centering
\includegraphics[width=0.6\textwidth]{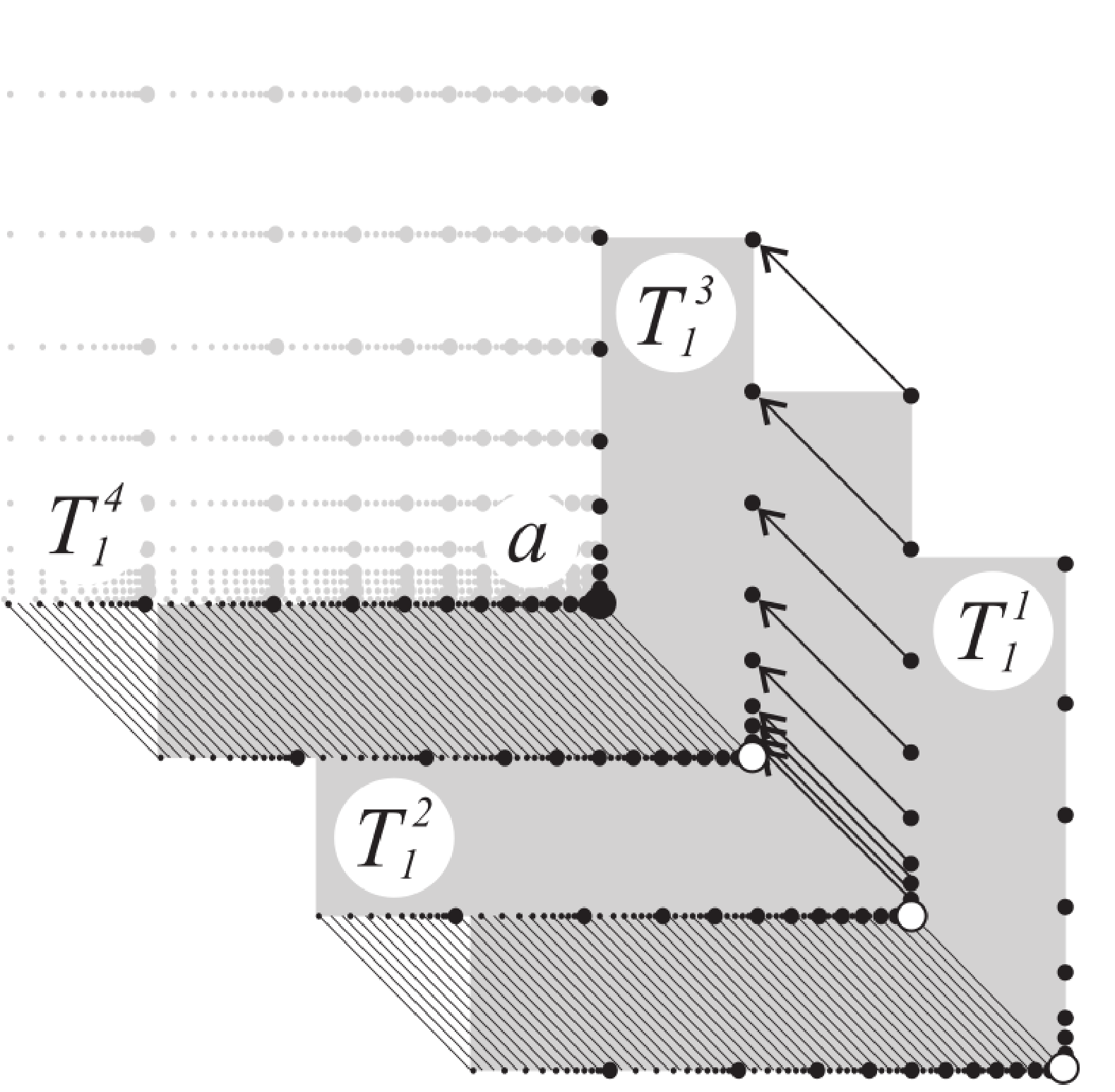} 
\caption{The space $G_4$}
\end{figure}




$\bullet$  $G_n$ is a $S(k)$-space for every $k\in \mathbb{N}$.

$\bullet$ $G_n$ is a $S(n)$-closed space. Indeed, consider any $S(n-1)$-cover
$V=\{V_{\alpha}\}$ of $G_n$, then the point
${\bf a}\in Int_{\theta^{n-1}} V_{\alpha}$ for some $V_{\alpha}$ from $V$.
Hence, $V_{\alpha}$ contains some $n$-hull of ${\bf a}$. But the complement of any closed $n$-hull of ${\bf a}$ is a compact set, so it is covered by a finite number of elements of $S(n-1)$-cover $V$.

On the other hand, the set $F=\{(1,\alpha,\omega_0)\}$ has a regular cardinality $\omega_1$ and does not have a point of
$\theta^0(n)$-complete accumulation ($G_n\setminus F$ is the
$n$-hull of the point ${\bf a}$ and it does not intersect $F$). Hence,

$\bullet$  $G_n$ is not weakly $S(n)$-$\theta$-closed.

So, $G_n$ is an example of an $S(n)$-space which
is $S(n)$-closed, but not weakly $S(n)$-$\theta$-closed.

\bigskip
Example \ref{126} completes the series of examples proving the strictness of all implications on Diagram 1. Moreover, Example \ref{124} and Example \ref{125} prove the independence of two classes of spaces in the $S(n)$-axiom of separation. Namely, Example \ref{124} proves that the class of weakly $S(n-1)$-closed spaces is not contained in the class of $S(n)$-$\theta$-closed spaces. And Example \ref{125}, on the contrary, proves that the class of $S(n)$-$\theta$-closed spaces is not contained in the class of weakly $S(n-1)$-closed spaces.

\bigskip

\begin{example}\label{noncom}

The unit interval $[0,1]$, endowed with the smallest topology which refines the euclidean topology, and contains $\mathbb{Q}\cap [0,1]$ as an open set is (weakly) $H$-closed but not compact.

\end{example}

 Thus, in the class of $S(n)$-spaces, the properties under consideration are in relations that are presented in the following diagram (Diagram 1):

\bigskip

\begin{figure}[h!]
\centering
\includegraphics[width=1\textwidth]{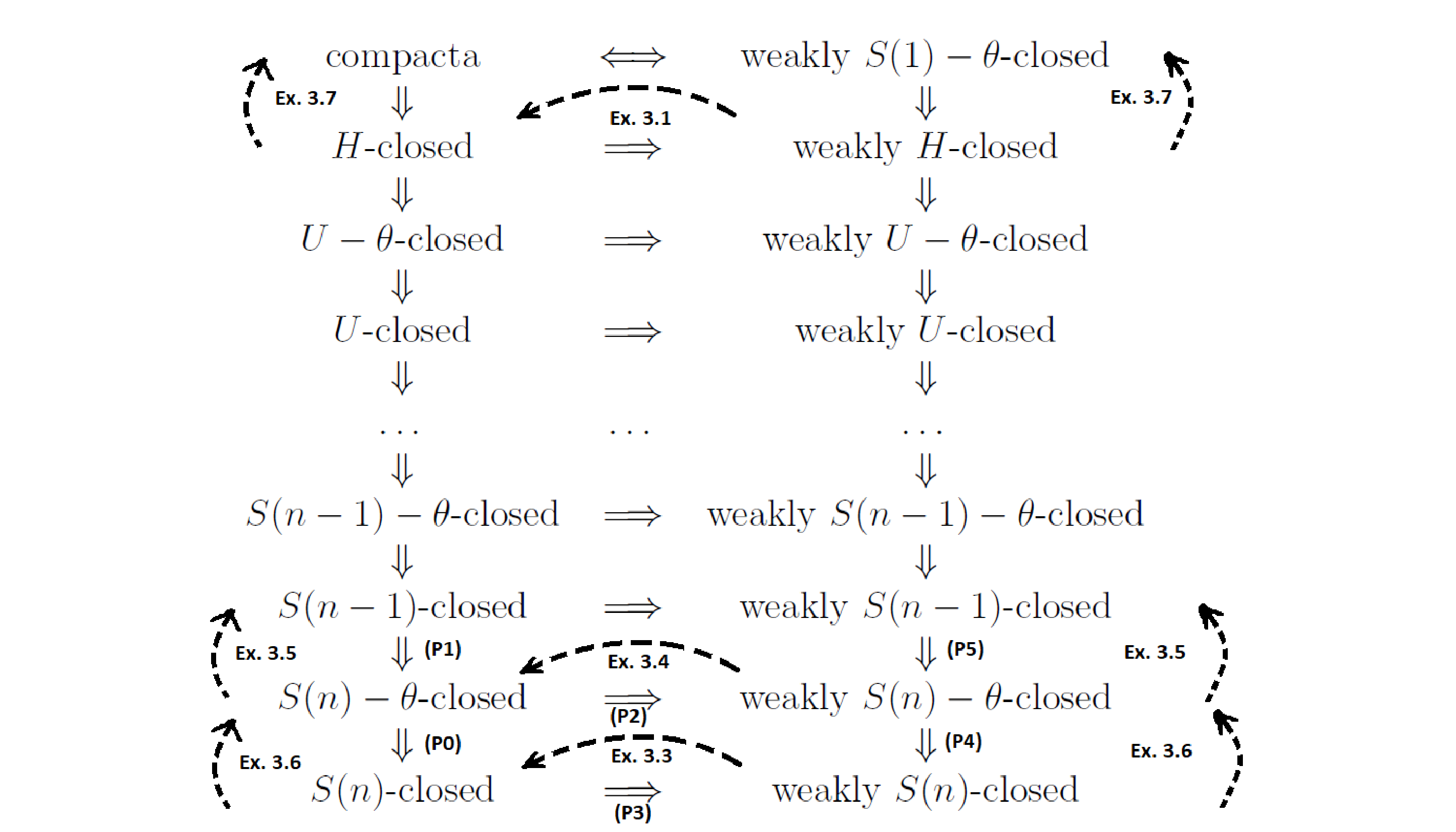}

Diagram 1.
\end{figure}




\newpage

\section{$S(n)$-closednesses in Lindel\"{o}f spaces}

Recall that a space $X$ is  {\it Lindel\"{o}f} if every open cover of $X$ has a countable subcover.

\bigskip

The following assertions are holds for every  $n>1$.

\bigskip

$({\bf L1})$  A Lindel\"{o}f weakly $S(n)$-closed $S(n)$-space is $S(n)$-closed (Theorem 3.6 in \cite{osip7}).

\bigskip
Recall that a topological space $X$ is called {\it linearly Lindel\"{o}f} if any uncountable (regular
cardinality) subset of $X$ has a point of complete accumulation. Note that every Lindel\"{o}f space is linearly Lindel\"{o}f.

$({\bf L2})$  A linearly Lindel\"{o}f weakly $S(n)$-$\theta$-closed $S(n)$-space is weakly $S(n-1)$-closed (Theorem 3.9 in \cite{osip7}).

\bigskip

$({\bf L3})$  A Lindel\"{o}f weakly $S(n)$-$\theta$-closed $S(n)$-space is $S(n-1)$-closed (Corollary 3.10 in \cite{osip7}).

\bigskip

$({\bf L4})$  A Lindel\"{o}f weakly $S(n)$-$\theta$-closed $S(n)$-space is $S(n)$-$\theta$-closed ({\bf L3} and {\bf P1}).

\begin{definition}(\cite{bcp}) Let $X$ be a space. For $n\in \mathbb{N}$, the $n$-$\theta$-closure of a subset $A$ of $X$ is $cl^{n}_{\theta}(A)=\underbrace{cl_{\theta}cl_{\theta}...cl_{\theta}(A)}_{n-times}.$

\end{definition}

\begin{definition}(\cite{bcp}) A space $X$ is a $\theta^n$-Urysohn, for every $n\in \mathbb{N},$ if for every $x,y\in X$ with $x\neq y$, there exist open subsets $U$ and $V$ of $X$ with $x\in U$ and $y\in V$ such that $cl^n_{\theta}(U)\cap cl^n_{\theta}(V)=\emptyset.$
\end{definition}

In (\cite{osip7}, Question 1) it posed the following question: {\it Does there exist a non
$S(n)$-$\theta$-closed Lindel\"{o}f $S(n)$-closed space for every
$n\geq 2$?}

The following example of non $\theta^n$-Urysohn $S(n)$-space answers of that question.

\begin{example} Fix $n\in \mathbb{N}$. Let
$\mathbb{R}=\bigcup\limits_{i=1}^{2n} A_i$ where $A_i$'s are
pairwise disjoint, each $A_i$ is dense in $\mathbb{R}$,
$|A_i|=\aleph_0$ for $i\neq 2$. Let $A'_{2n+1}$ be a copy of $A_1$
and let $X_{2n+1}=\bigoplus\limits_{i=1}^{2n+1} A_i$.

If $a,b\in \mathbb{R}$ and $a<b$, an open base for $X_{2n+1}$ is
generated be the following families of sets:

(1)  $(a,b)\cap A_i$ for all odd $1\leq i\leq 2n+1$;

(2)  $(a,b)\cap (A_{i-1}\cup A_i\cup A_{i+1})$ for all even $2\leq
i\leq 2n$.

\end{example}

\begin{figure}[h!]
\centering
\includegraphics[width=0.8\textwidth]{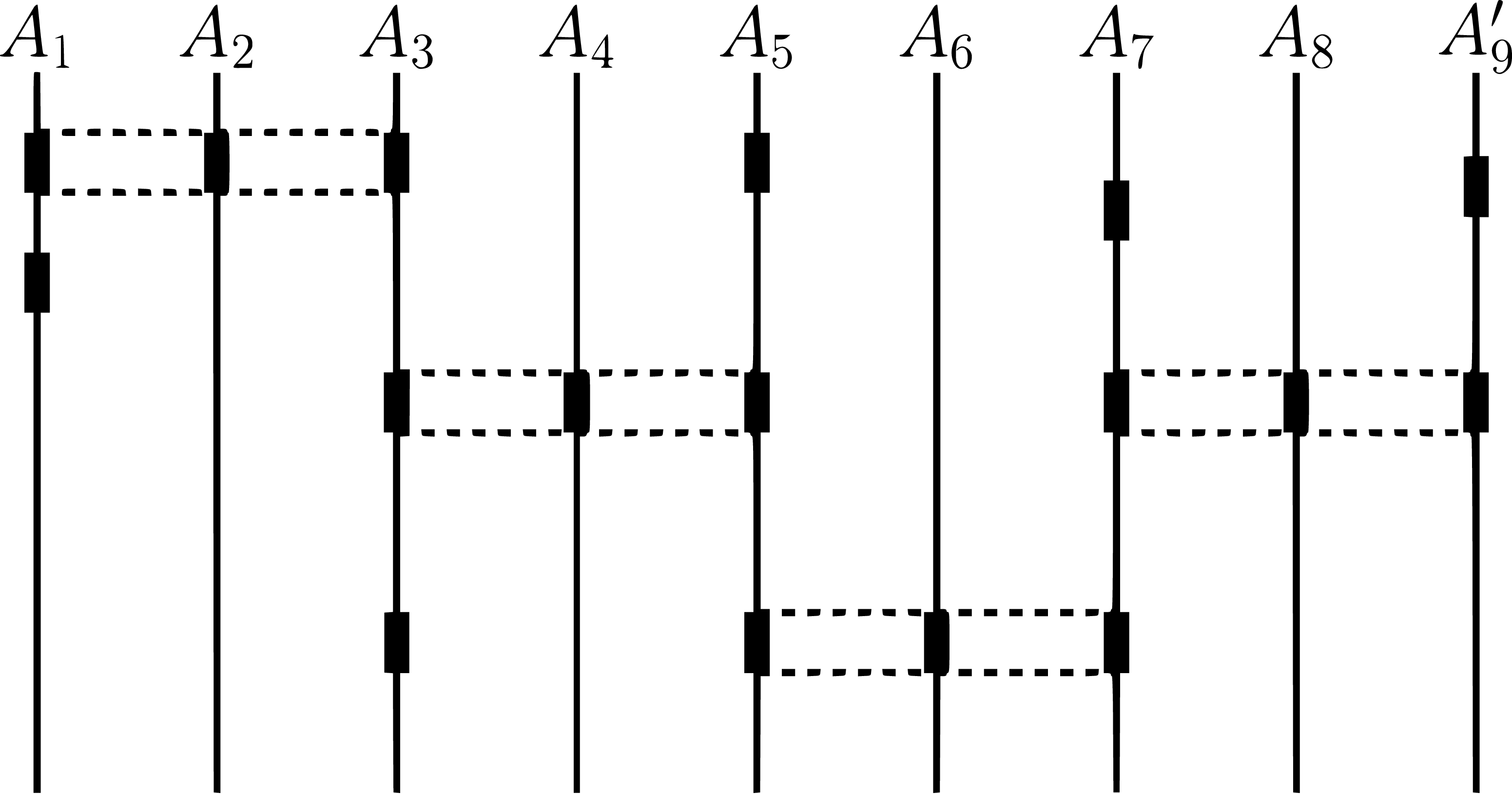}

\caption{The space $P_4$}
\end{figure}



The obtained space $P_n$ (see Figure~7 for $n=4$) has the following
properties.

Let $a,b,c,d\in \mathbb{R}$, $U=(a,b)\cap A_1$. Then $cl_{\theta}(U)=cl(U)\subseteq [a,b]\cap (A_1\cup A_2).$  Let $V=(c,d)\cap (A_3\cup A_4\cup A_5).$ Then $cl_{\theta}(V)=cl(V)=[c,d]\cap (A_2\cup A_3\cup A_4\cup A_5\cup A_6)$. It follows that $cl^2_{\theta}(U)=[a,b]\cap (A_1\cup A_2\cup A_3\cup A_4)$. By induction, $cl^n_{\theta}(U)=[a,b]\cap (A_1\cup A_2\cup ... \cup A_{2n})$. Likewise, starting from the right-hand subspace $A'_{2n+1}$ with $U=(a,b)'\cap A'_{2n+1}$, we have $cl^n_{\theta}(U)=[a,b]'\cap(A'_{2n+1}\cup A_{2n}\cup ... \cup A_2)$.

We have the following consequences:

$\bullet$  $P_n$ is a $S(n)$-space. Every $x,y\in P_n$, $x\neq y$, are $\theta^n$-separated, i.e., there are $n$-hull $V(x)$ of $x$ and $n$-hull $W(y)$ of $y$ such that
if $n$ is odd then $V(x)\cap W(y)=\emptyset$, and
if $n$ is even then $\overline{V(x)}\cap \overline{W(y})=\emptyset$.

$\bullet$ Every pair of points $x$, $x'$ such that $x\in A_1$, $x'\in A'_{2n+1}$ ($x'$ is a
copy of $x$) are not $\theta^n$-Urysohn separated. Let $a$, $b$, $c$, $d\in \mathbb{R}$, $x\in (a, b)$ and
$x'\in (c, d)'$, $U = (a, b)\cap (c, d)$. Then $cl^n_{\theta}(U)=\overline{U}\cap (A_1\cup A_2\cup ... \cup A_{2n})$ and

$cl^n_{\theta}(U')=\overline{U'}\cap (A'_{2n+1}\cup A_{2n}\cup ... \cup A_{2})$. Thus, $cl^n_{\theta}(U)\cap cl^n_{\theta}(U')\neq \emptyset$; $P_n$ is not $\theta^n$-Urysohn.

Consider the subspace $S=[0,1]\bigcap (\bigcup\limits_{i=1}^{2n} A_i)$ of $P_n$.

Then $S$ is a Lindel\"{o}f $S(n)$-closed space, but it is not $S(n)$-$\theta$-closed space.

1. Since $[0,1]\cap A_2$ is subspace of $\mathbb{R}$ and $\mathbb{R}$ is hereditarily Lindel\"{o}f,
$[0,1]\cap A_2$ is Lindel\"{o}f and, hence, $S$ is Lindel\"{o}f.

2. Let $a\in [0,1]\cap A_1$. Consider a sequence $\{a_m : m\in \mathbb{N}\}$ such that
$a_m \in [0,1]\cap A_{2n}$ for every $m\in \mathbb{N}$ and $\{a_m\}_{m\in \mathbb{N}}$ converges to $a$ ($m\rightarrow \infty$)
in natural topology of $[0, 1]$. Then there is a $n$-hull $U(a)$ of the point $a$ such
that $U(a)\bigcap \{a_n : n\in \mathbb{N}\}=\emptyset$. It follows that the set $\{a_n : n\in \mathbb{N}\}$ has not
a $\theta^0(n)$-complete accumulation point. Hence, $S$ is not $S(n)$-$\theta$-closed space.

3. Note that $S$ is weakly $S(n)$-closed space. Then, by ${\bf (L1)}$, $S$ is
a $S(n)$-closed space.

\bigskip

Thus,  classes of the considered spaces in Lindel\"{o}f $S(n)$-spaces
are presented in the following diagram (Diagram 2):

\bigskip

\begin{figure}[h!]
\centering
\includegraphics[width=0.7\textwidth]{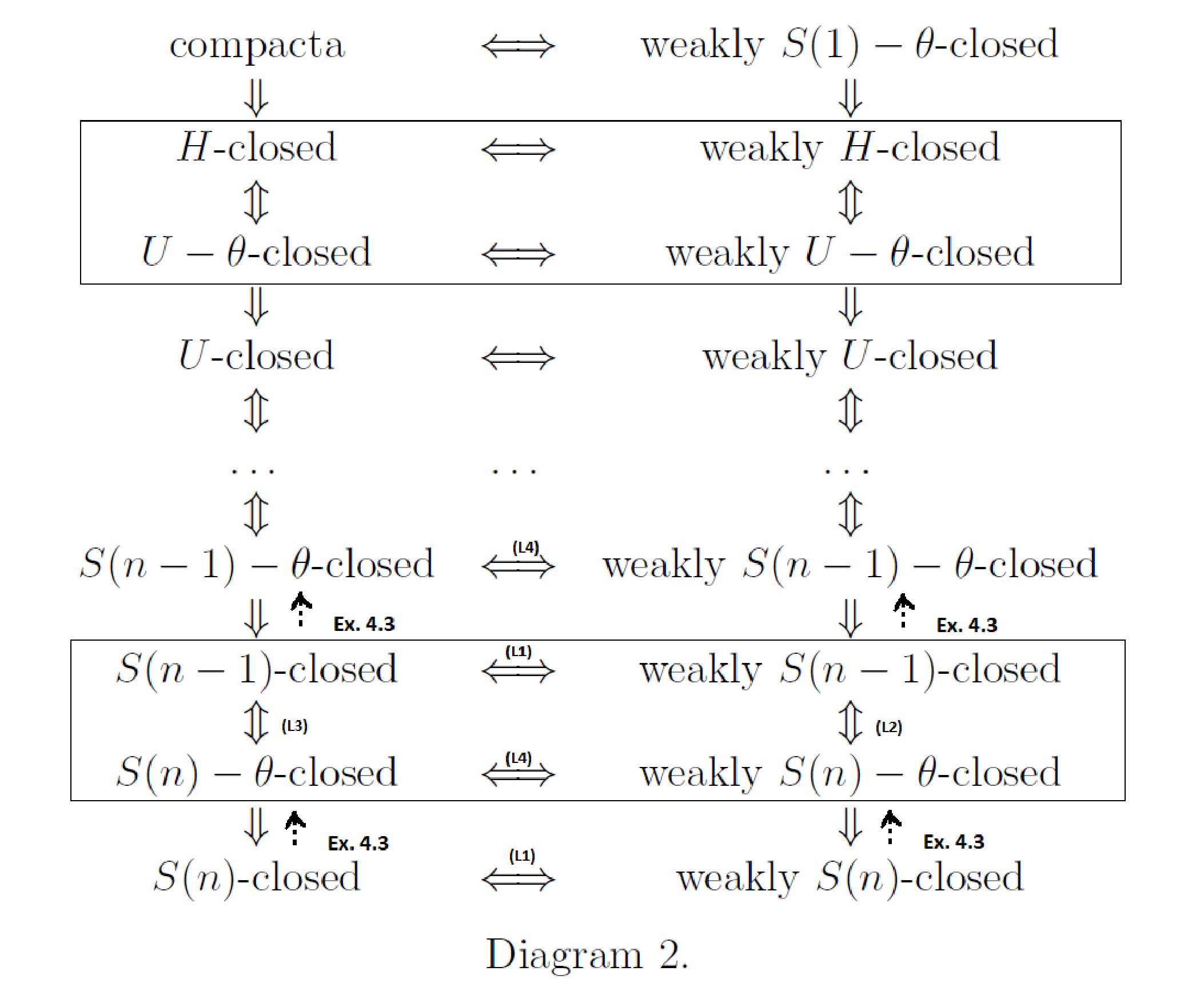}

\end{figure}


\section{Solving some open problems with examples}

This section will present solutions to some of the problems posed in \cite{dg}, \cite{bps}. These solutions were announced in the papers  \cite{osip4,osip16}.

\subsection{Question 1}

By $\tau_{\theta}$ we denote the topology obtained from $\tau$ by
declaring all $\theta$-closed sets in $\tau$ to be closed.

Consider the following question (Question 2.6 in  \cite{dg}):

\medskip

{\it Does there exist a Urysohn-closed space $(X,\tau)$ with $(X,\tau_{\theta})$ quasi-compact such that $(X,\tau)$ is not Urysohn-$\theta$-closed?}

\medskip

The following example answers of that question.

\begin{example}\label{311}  Let $\omega_1$ be the set of all countable ordinals, and let $\omega_1+1=\omega_1\cup\{\omega_1\}$. We denote the points of $\omega_1+1$ by $\{\alpha\}$ and the points
of $\omega_1$ by $\{\alpha^1\}$. Let $P$ be the set of all limit points of $\omega_1$ in the order topology, and let $P_1$  be the set of all limit points of  $P$. We strengthen the order
topology at the points of  $P_1$  as follows. For the base neighborhoods of each point $\alpha^1\in P_1$  we take
  $U(\alpha^1)=(O(\alpha^1)\setminus P)\cup \{\alpha^1\}$, where
  $O(\alpha^1)$ is a neighborhood in the order topology. In the topological sum $(\omega_1+1)\bigoplus\omega_1$, we identify $\alpha$ with $\alpha^1$ for each $\alpha^1\in P_1$ and denote
   the points glued together by $\{\alpha, \alpha^1\}$.

   In the quotient space thus obtained, we strengthen the topology at the point  $\{\omega_1\}$. For its base
neighborhoods we take

 $U(\omega_1)=(O(\omega_1)\setminus P_1^1)$, where $O(\omega_1)$ is any neighborhood of $\{\omega_1\}$ in the order
topology and $P_1^1$ is the image of  $P_1$  under the quotient map.

\bigskip
We denote the resulting quotient space by $F_1$ (see Figure~8 where ${\bf a}$  is the point $\{\omega_1\}$).

\end{example}

\begin{figure}[h!]
\centering
\includegraphics[width=0.8\textwidth]{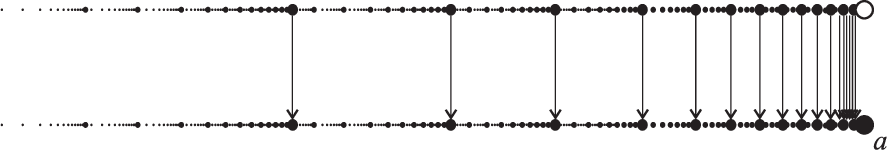}

\caption{The space $F_1$}
\end{figure}



\bigskip

The space $F_1$ has the following properties:

$\bullet$ $F_1$ is an $S(n)$-space for any $n\in \mathbb{N}$;

This follows from the fact that for any two points $x$ and $y$, there exists a $2$-hull $O_x$ of $x$ such that
$y\notin O_x$ and $O_x$ is an open-closed subset of the space $F_1$.

$\bullet$ $F_1$ is $S(2)$-closed ($U$-closed);

Indeed, take any $S(1)$-covering $U=\{U_s\}$ of the space $F_1$.
Then the point $\{\omega_1\}$ lies in the $\theta$-interior
of some $U_s$ from $U$. Thus, $U_s$ is the
$2$-hull of $\{\omega_1\}$, but for
$F_1\setminus\overline{U_s}$, there obviously exists a finite set
$U_{s_1},...,U_{s_k}$ such that
$F_1\setminus\overline{U_s}\subset\bigcup\limits_{i=1}^k\overline{U_{s_i}}$.
 Thus, the set $\overline{U_s}$,$\overline{U_{s_1}}$,...,$\overline{U_{s_k}}$ of $U$
covers the space $F_1$.

$\bullet$ $F_1$ is not an $S(2)$-$\theta$-closed (not $U$-$\theta$-closed) space.

 Since the set $P\setminus P_1$, having regular cardinality $\omega_1$, does not have a
$\theta^0(2)$-accumulation point, the space $F_1$ is not a weakly
$S(2)$-$\theta$-closed space and (by $({\bf P2})$) is not
$S(2)$-$\theta$-closed.

$\bullet$ $(F_1,\tau_{\theta})$ is compact.

 Indeed, consider the $\tau_{\theta}$ topology on $F_1$.
 Let $V$ be an arbitrary $\tau_{\theta}$-neighborhood of the point
 $\{\omega_1\}$. Then there exists $\alpha<\omega_1$ such that
 $\{\beta:\beta>\alpha\}\subset V$. Consequently, $V$
 also contains the glued points $(\beta,\beta^1)$ where $\beta^1\in
 P_1$, $\beta>\alpha$. For each such point $(\beta,\beta^1)$
 there exists a neighborhood $V_{(\beta,\beta^1)}\subset V$ due to the
 $\theta$-closure of the set $F_1\setminus V$. This means that
 the set $\{\beta^1:\beta^1>\alpha, \beta^1<\omega_1\}$
 also lies in $V$. It easily follows from this that the space
 $(F_1,\tau_{\theta})$ is compact.

\bigskip

Thus,  the $U$-closed space $F_1$ is not (weakly) $U$-$\theta$-closed, but $(F_1, \tau_{\theta})$ is compact.

\subsection{Question 2}

An open set $U$ of a topological space $(X,\tau)$ is {\it regularly open} if $U=int \overline{U}$. Recall that a topology is said to be {\it semiregular} if it has a base consisting of regular open sets.
By $\tau_s$ we denote the topology obtained from $\tau$ by declaring a set to be closed in $\tau_s$ if and only if it
is regular closed in $\tau$.

\medskip

It is easy to see that a space $(X, \tau)$ is $S(n)$-closed whenever its semiregularization $(X, \tau_s)$ is $S(n)$-closed. On the other hand there exist non compact $H$-closed Urysohn spaces (see Example \ref{noncom}) (their semiregularization is always compact).

\bigskip
Consider the following question (Problem 2 in \cite{dg}):

\medskip

{\it Is it true that any
$S(n)$-space $(X,\tau)$ for which $(X, \tau_s)$ is $S(n)$-$\theta$-closed is $S(n)$-$\theta$-closed?}

\medskip

Example \ref{311} (independely of \cite{jrw}) answers of that question for $n=1$.

\medskip

The following example answers of that question for any $n>1$.

\begin{example}\label{314} Let $n>1$ and $T^i$ ( $i=1,...,n-1$ ) be  $n-1$
copies of the deleted Tychonoff plane $T$, whose elements will be denoted
$(i,\alpha,k)$. Let $\omega_1$ be the set of all countable
ordinals with the order topology, whose elements will be denoted
by  $\{\alpha\}$. Let $P$ be the set of all limit points of $\omega_1$ in the order topology, and let $P_1$  be the set of all limit points of  $P$. We strengthen the order
topology at the points of  $P_1$  as follows.

 For the base neighborhoods of each point $\alpha\in P_1$  we take
  $U(\alpha)=(O(\alpha)\setminus P)\cup \{\alpha\}$, where
  $O(\alpha)$ is a neighborhood in the order topology.

Let $\varphi: \omega_1\rightarrow P_1$ be a bijection (order-preserving) function.

 Consider the following identifications on the topological sum $\bigoplus\limits_{i=1}^{n-1}T^i\bigoplus\omega_1$:

  $(1,\alpha,\omega_0)\sim (\varphi(\alpha))$ for every  $\alpha\in \omega_1$;

 $(s,\omega_1,k)\sim (s+1,\omega_1,k)$ for all odd $s$;

 $(s,\alpha,\omega_0)\sim (s+1,\alpha,\omega_0)$ for all even  $s$.

To the obtained space, we add the point ${\bf a}$ with the base of
neighborhoods:

 $U_{\alpha,k}({\bf a})=\{(n-1,\beta,p): \alpha < \beta <\omega_1, k<p<\omega_0
 \}\cup\{{\bf a}\}$.

We denote the resulting quotient space by $F_n$ (see Figure~9 for $n=3$).

\end{example}

\begin{figure}[h!]
\centering
\includegraphics[width=0.7\textwidth]{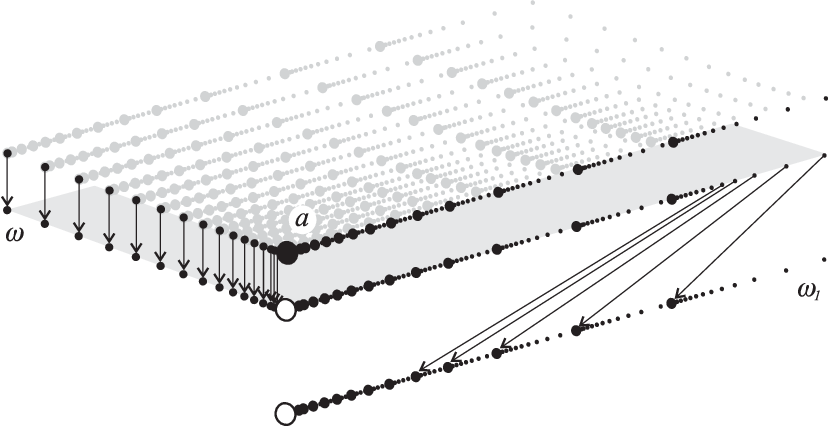}

\caption{The space $F_3$}
\end{figure}




The space $F_n$ has the following properties:

$\bullet$   $F_n$ is a $S(k)$-space for every $k\in \mathbb{N}$.

$\bullet$  $F_n$ is a $S(n)$-closed space.

On the other hand, the set $F=P\setminus P_1$ in $\omega_1$
has regular cardinality $\omega_1$ and does not have a
$\theta^0(n)$-accumulation point. Therefore,

$\bullet$ $F_n$ is not weakly $S(n)$-$\theta$-closed space
and, by ${\bf (P2)}$,  is not $S(n)$-$\theta$-closed.

$\bullet$ $(F_n,\tau_s)$ is a $S(n)$-$\theta$-closed space.

 Indeed, $\omega_1$ in $(F_n,\tau_s)$ has the
 order topology. Therefore, any neighborhood of the set
 $Z=P_1\setminus S$, where $S$ is countable, will contain all points of $\omega_1$ but a countable number. Note that in $(F_n,\tau_s)$ the point ${\bf a}$ will be a
 $\theta^0(n)$-complete accumulation point for any uncountable subset
$\omega_1$. Consider an arbitrary $n$-hull $U_a$ in $(F_n,\tau_s)$
 of the point ${\bf a}$, then $F_n\setminus U_a$ is compact. Therefore,
 $(F_n,\tau_s)$ is an $S(n)$-$\theta$-closed space.

Thus, $F_n$ is an example of an $S(n)$-space that
is $S(n)$-closed, but not $S(n)$-$\theta$-closed
space, which in the semiregular topology $\tau_s$
is an $S(n)$-$\theta$-closed space.

\bigskip

Note that Example \ref{314} also solves negatively the question of
relatively weakly $S(n)$-$\theta$-closedness. Indeed,
$(F_n,\tau_s)$ is weakly $S(n)$-$\theta$-closed, but
 the space $(F_n,\tau)$ is not a weakly $S(n)$-$\theta$-closed space.

\subsection{Question 3}

A regular space $X$ is called regular-closed if it is a closed subspace in every regular
space in which it is embedded. A topological space $X$ is feebly compact if any open locally finite
family of its subsets is finite.

\medskip

In 1982, Pettey \cite{pett} proved that the product of regular-closed spaces is regular-closed if it is
feebly compact. The validity of a similar theorem in the class of $U$-$\theta$-closed spaces was discussed
in \cite{dg}, where the problem (Problem 5) on the product of $U$-$\theta$-closed
spaces  was formulated:

\medskip

 {\it It is required to prove or to disprove that the product
of $U$-$\theta$-closed spaces is feebly compact.}

\medskip
Note that a Lindel\"{o}f $U$-$\theta$-closed space is $H$-closed (see ${\bf (L3)}$). It was by Chevalley and Frink that products of $H$-closed spaces are $H$-closed. Observe that every $H$-closed space is feebly compact.
It follow that for the case of Lindel\"{o}f $U$-$\theta$-closed
spaces that problem is solved positively.

\medskip
Note that $U$-$\theta$-closedness is not a multiplicative property
~\cite{dg}.

\medskip

 Next, two Urysohn  $U$-$\theta$-closed spaces are constructed such that their product is not
feebly compact, and thus the general question of the feebly compactness
of the product of two $U$-$\theta$-closed spaces is answered
negatively.

\medskip

For the construction, we use the construction of example 3.10.19 from
the book~\cite{eng}.

\begin{example} Let $\beta \mathbb{N}$ be the Stone-\v{C}ech
extension of the set of positive integers $\mathbb{N}$. For every
$M\subset \beta \mathbb{N}$, denote by $\mathcal{P}(M)$ the family of all
countable infinite subsets of the set $M$; let $f$ be the mapping
which to every term $S$ of the family $\mathcal{P}(\beta \mathbb{N})$ puts
into correspondence some limit point of the set $S$ in the space
$\beta \mathbb{N}.$

 Setting $X_0=\mathbb{N}$ and $X_\alpha=\bigcup\limits_{\gamma < \alpha} X_\gamma\bigcup
f(P(\bigcup\limits_{\gamma < \alpha} X_\gamma))$ for
$0<\alpha<\omega_1$,  we define the transfinite sequence $
X_0$, $X_1$, ..., $X_{\alpha},...,$ $\alpha<\omega_1$, of subsets of the
space  $\beta \mathbb{N}$ by means of the transfinite induction.
The space $A=\bigcup\limits_{\alpha<\omega_1} X_\alpha$ is
countably compact \cite{eng} because every $S\in \mathcal{P}(A)$ is contained
in some $X_{\alpha}$, and consequently, it has a limit point in
$X_{\alpha+1}$ and all the more in $A$. It is easily shown by the
transfinite induction that $|X_\alpha|\le \mathfrak{c}$. Hence,
 $|A|\le \mathfrak{c}$.

We set $B=\mathbb{N}\cup(\beta \mathbb{N}\setminus A)$. Since the
power of the closure of $S$ equals $2^{\mathfrak{c}}$ for every
$S\in \mathcal{P}(B)$, each infinite subspace of the space $B$ has a limit
point in $B$, and the space $B$ is countably compact.

Thus, we have the partition of the remainder $\mathbb{N}^*=\beta
\mathbb{N}\setminus \mathbb{N}$ into two disjoint subsets
$A_1=A\setminus \mathbb{N}$ and $B_1=B\setminus \mathbb{N}$ such
that $A$ and $B$ are countably compact.

 Let $X=(\beta \mathbb{N}\times [0,1])\setminus(\mathbb{N}\times (0,1])$.

\end{example}

Note that $X$ (with the topology of product) is a compact space
(as closed subspace of the product $\beta \mathbb{N}\times
[0,1])$. We define a stronger topology of $X$:

If $x\in A_1$, then the base of the point $(x,0)$ is the family of
all sets of the form
      $U(\{x,0\})=\{x,0\}\bigcup O(\{x,0\})\setminus(\beta \mathbb{N}\times \{0\})$,
where $O(\{x,0\})$ is an arbitrary neighborhood (in the topology
of product) of the point $(x,0)$ in $X$. The bases of other points
are the same as in the topology of product. We denote by $D_1$ the
topological space obtained.

Let us define one more topology of $X$:

If $x\in B_1$, then the base of the point $(x,0)$ is the family of all sets of the form $U(\{x,0\})=\{x,0\}\bigcup O(\{x,0\})\setminus (\beta \mathbb{N}\times \{0\})$, where $O(\{x,0\})$ is an arbitrary neighborhood (in the topology of product) of the point $(x, 0)$ in $X$. The bases of other points are the same as in the topology of
product. We denote by $D_2$ the topological space obtained.

\begin{figure}[h!]
\centering
\includegraphics[width=0.8\textwidth]{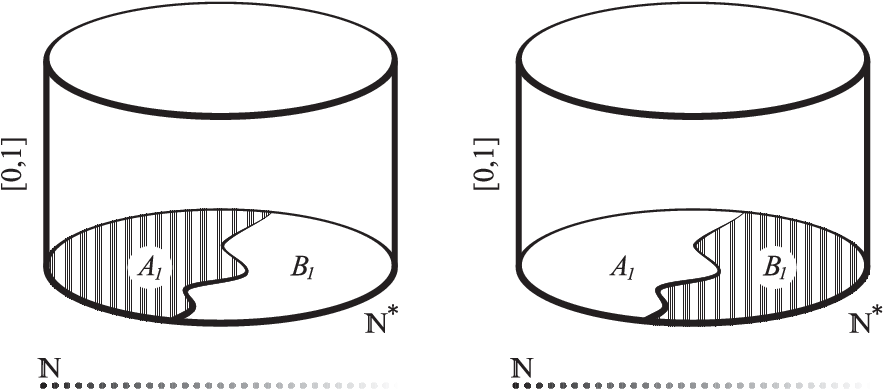}

\caption{The spaces $D_1$ and $D_2$}
\end{figure}



The spaces $D_1$ and $D_2$ are Urysohn, since they are condensed on the Hausdorff compact space $X$.
Let us prove that $D_1$ and $D_2$ are $U$-$\theta$-closed spaces. Indeed, let $\{V_{\alpha}\}_{\alpha}$ be an Urysohn cover, for
example, of the space $D_1$. Then, for any point $x\in S$  $(S=\mathbb{N}^*\times [0,1])$, there exists a neighborhood
$O_x$ which is open in the topology of product and such that $O_x\subset V_{\alpha}\bigcap S$ for some $\alpha$. Then,
$\{O_x\}_{x\in S}$ covers the compact set $S$; hence, there exists a finite subcover $O_{x_1}$, $O_{x_2}$, ..., $O_{x_k}$ and, as
a consequence, a finite subcover $V_1\cap S$, $V_2\cap S$, ..., $V_k\cap S$ of the space $S$. Since $(B_1\times \{0\})\subset S$, and $V_1\cap S$, $V_2\cap S$, ..., $V_k\cap S$ cover of $B_1\times \{0\}$, we get that
$V_1$, $V_2$, ..., $V_k$ cover of $\mathbb{N}\times \{0\}$, which
that $V_1$, $V_2$, ..., $V_k$ is a finite subcover $D_1$. The $U$-$\theta$-closedness of $D_2$ is proved
analogously.

Consider the Tychonoff product of the two $U$-$\theta$-closed spaces $D_1$ and $D_2$ : $Y=D_1\times D_2.$ Denote
by $\Delta$ the diagonal of the product and by $\Delta_0$ a subset of $\Delta$, $\Delta_0=\{((1,0), (1,0)), ((2,0), (2,0)), ...\}$.
Since $D_1$ and $D_2$ are Urysohn, the diagonal $\Delta$ is closed in $Y$. Let $a=\{(x,y),(x,y)\}\in \Delta\setminus \Delta_0$. If
$y\neq 0$, then, obviously, the point $a$ has a neighborhood $V_a$ such that $V_a\cap \Delta_0=\emptyset$.
Let us show that such a neighborhood exists also in the case $y=0$. Indeed, let $a=((x,0), (x, 0))\in \Delta\setminus \Delta_0;$ then,
$x\in A_1$ or $x\in B_1$. If $x\in A_1$, then in $D_1$ there exists (by construction) a neighborhood $U(\{x, 0\})$
such that $U(\{x,0\})\cap (\mathbb{N}\times \{0\})=\emptyset$ and $U(\{x,0\})\times O(\{x,0\})\cap \Delta_0=\emptyset$ for any neighborhood
$O(\{x, 0\})$  in $D_2$.

If $x\in B_1$, then there exists (by construction) a neighborhood $U(\{x,0\})$ in $D_2$ such that
$U(\{x, 0\})\cap (\mathbb{N}\times \{0\})=\emptyset$ and $U(\{x,0\})\times O(\{x,0\})\cap \Delta_0=\emptyset$ for any neighborhood
$O(\{x, 0\})$  in $D_1$.

Thus, the point $a$ has the desired neighborhood $V_a$. So, $\Delta_0$ is closed in $\Delta$ and, hence, is closed in
the product $Y$.

The set $\{(i,0)\}$ is open in $D_1$ and in $D_2$ for $i=1, 2,...;$, hence, $\Delta_0$ is an open subset in $Y$. The
locally finite open system $\{(i,0),(i,0)\}_{i\in \mathbb{N}}$ is not finite, and, therefore, $Y$ is not feebly compact.
Note that the space $Y$ is not countably compact because it contains the infinite discrete closed set $\Delta_0$.

\bigskip

 Since $\Delta_0$ is an infinite discrete closed set,
which does not have a point $x$ such that $|\Delta_0\bigcap V_x |=|
\Delta_0 |$ for any 2-hull of $x$, $Y$ is not
weakly $U$-$\theta$-closed and, hence, it is not
$U$-$\theta$-closed.

\bigskip

 Thus, the constructed space $Y$ proves that the
property of weakly $U$-$\theta$-closedness is not multiplicative.
We have obtained that the problem of feebly compactness of the product of $U$-$\theta$-closed
spaces, in general, is solved negatively.

\bibliographystyle{model1a-num-names}
\bibliography{<your-bib-database>}

\end{document}